# SEPARATION CUT-OFFS FOR BIRTH AND DEATH CHAINS

By Persi Diaconis[1] and Laurent Saloff-Coste[2]

*Stanford University and Cornell University*

This paper gives a necessary and sufficient condition for a sequence of birth and death chains to converge abruptly to stationarity, that is, to present a cut-off. The condition involves the notions of spectral gap and mixing time. Y. Peres has observed that for many families of Markov chains, there is a cut-off if and only if the product of spectral gap and mixing time tends to infinity. We establish this for arbitrary birth and death chains in continuous time when the convergence is measured in separation and the chains all start at 0.

**1. Introduction.** Some ergodic Markov chains show a sharp transition in convergence to stationarity. This was first observed for random transpositions on the symmetric group in [22]. The phenomenon was clearly identified in [2] where the term "cut-off phenomenon" was introduced (see [2], Figure 2). Recently, Yuval Peres observed that for many examples, a cut-off occurs if and only if the product $\lambda\tau$ tends to infinity where $\lambda$ is the spectral gap (i.e., 1 minus the second largest eigenvalue) and $\tau$ is the mixing time (i.e., the first time the distance to stationarity is less than $1/4$).

Our main theorem proves a precise version of this statement for all finite continuous-time birth and death chains started at 0 when convergence is measured in separation distance. Namely, for each $n$, let $\gamma_n^t$ denote the distribution at time $t > 0$ of a given ergodic birth and death chain on $\Omega_n = \{0, 1, \ldots, n\}$, started at 0. Let $\nu_n$ be the corresponding stationary measure. The separation between $\gamma_n^t$ and its target $\nu_n$ is

$$\mathrm{sep}(\gamma_n^t, \nu_n) = \sup_{x \in \Omega_n} \{1 - \gamma_n^t(x)/\nu_n(x)\}.$$

Received June 2003; revised March 2006.
[1]Supported in part by NSF Grant DMS-05-05673.
[2]Supported in part by NSF Grant DMS-01-02126.

Let $\lambda_{n,i} \in [0,2]$, $i = 0, \ldots, n$, be the eigenvalues, in nondecreasing order, of minus the associated $Q$-matrix $Q_n$ (i.e., the matrix representing the infinitesimal generator of the associated Markov process). We always have $\lambda_{n,0} = 0$, and $\lambda_n = \lambda_{n,1}$ is called the *spectral gap*. Set

$$t_n = \sum_1^n \lambda_{n,i}^{-1}.$$

THEOREM 1.1. *Referring to the setting and notation introduced above, we have:*

1. *If $\lambda_n t_n$ tends to infinity with $n$, then, for any $c \in (0,1)$,*

$$\lim_{n \to \infty} \mathrm{sep}(\gamma_n^{(1+c)t_n}, \nu_n) = 0, \qquad \lim_{n \to \infty} \mathrm{sep}(\gamma_n^{(1-c)t_n}, \nu_n) = 1.$$

2. *If $\lambda_n t_n$ does not tend to infinity with $n$, then, for any $c \in (0,1)$,*

$$\limsup_{n \to \infty} \mathrm{sep}(\gamma_n^{(1+c)t_n}, \nu_n) > 0, \qquad \liminf_{n \to \infty} \mathrm{sep}(\gamma_n^{(1-c)t_n}, \nu_n) < 1.$$

See Theorem 5.1 below. A detailed analysis of the cut-off window is also obtained. The proof uses the duality theory of [16] to convert convergence rates into first hitting time estimates and Keilson's representation of first hitting times as sums of independent exponentials with parameters related to the spectrum of the chain.

The paper is organized as follows. Section 2 discusses various distances and carefully defines the cut-off phenomenon. Section 3 gathers elementary remarks concerning the cut-off phenomenon in separation and total variation distance. Birth and death chains are introduced in Section 4 which reviews duality and Keilson's spectral representation of hitting times. The main results—Theorems 5.1 and 5.2—are stated and proved in Section 5. They provide a characterization of the cut-off phenomenon for continuous-time birth and death chains started at 0 (in discrete time, the result is restricted to chains satisfying a certain monotonicity condition). Section 6 gives, when it exists, a precise description of the shape of the separation cut-off. This shape may or may not be Gaussian. It is Gaussian if and only if the size of the window is of an order of magnitude strictly larger than the relaxation time $1/\lambda$ (i.e., the inverse of the spectral gap). Section 7 gives detailed examples comparing cut-offs in separation, total variation and $L^2$-distance. These examples includes simple random walks, Metropolis chains, the Bernoulli–Laplace and Erhrenfest chains and simple random walk on distance-transitive graphs.

**2. Distances and cut-offs.** Assume that to any finite set $\Omega$ and any pair of probability measures $\mu, \nu$ on $\Omega$ is associated a real number $D(\mu, \nu)$ such that $D(\mu, \nu) \in [0, 1]$,

$$\sup_{\Omega} \sup_{\mu, \nu} D(\mu, \nu) = 1$$

and $D(\mu, \nu) = 0$ if and only if $\mu = \nu$ (in some cases, the supremum in $\Omega$ in the equation displayed above might be necessary to attain the supremum of the possible values taken by $D$). Examples of interest are the total variation distance

$$D(\mu, \nu) = \|\mu - \nu\|_{\mathrm{TV}} = \sup_{A \subset \Omega} \mu(A) - \nu(A)$$

and separation

$$D(\mu, \nu) = \mathrm{sep}(\mu, \nu) = \max_{\omega \in \Omega} \left\{ 1 - \frac{\mu(\omega)}{\nu(\omega)} \right\}.$$

Note that separation is not symmetric and is not a distance between probability measures. Separation was introduced in [2, 3] in the context of the study of convergence of ergodic finite Markov chains.

Consider a sequence of (finite) probability spaces $(\Omega_n, \nu_n)$, $n = 1, 2, \ldots$, each equipped with a sequence of probability measures $\mu_n^k$, $k = 0, 1, 2, \ldots$, such that

$$\lim_{k \to \infty} D(\mu_n^k, \nu_n) = 0.$$

DEFINITION 2.1. A family $(\Omega_n, \nu_n, (\mu_n^k)_{k=0,1,\ldots})_{n=1,2,\ldots}$ presents a cut-off (more precisely, a $D$-cut-off) if there exists a sequence $(t_n)$ of positive reals such that for any $\varepsilon \in (0, 1)$,

(a) $\lim_{n \to \infty} D(\mu_n^{k_n}, \nu_n) = 0$ if $k_n > (1 + \varepsilon) t_n$ for all sufficiently large $n$,
(b) $\lim_{n \to \infty} D(\mu_n^{k_n}, \nu_n) = 1$ if $k_n < (1 - \varepsilon) t_n$ for all sufficiently large $n$.

The next definition introduces the notion of window size for the cut-off phenomenon:

DEFINITION 2.2. Given sequences $(t_n)$ and $(b_n)$ of positive reals, we say that the family

$$(\Omega_n, \nu_n, (\mu_n^k)_{k=0,1,\ldots})_{n=1,2,\ldots}$$

presents a $(t_n, b_n)$-cut-off [more precisely, a $(t_n, b_n)$-$D$-cut-off] if $b_n/t_n$ tends to zero and

(a) $f_+(c) = \limsup_{n \to \infty} D(\mu_n^{\lceil t_n + c b_n \rceil}, \nu_n)$ satisfies $\lim_{c \to \infty} f_+(c) = 0$,

(b) $f_-(c) = \liminf_{c \to \infty} D(\mu_n^{\lfloor t_n - cb_n \rfloor}, \nu_n)$ satisfies $\lim_{c \to \infty} f_-(c) = 1$.

Both definitions can be interpreted in an obvious way when the discrete family of measure $\mu_n^k$, $k = 0, 1, \ldots$, is replaced by a continuous family $\mu_n^t$, $t \geq 0$ (in this case, $f_+$ and $f_-$ are defined using $t_n + cb_n$ without rounding to the next or previous integer). Versions of these definitions were introduced in [3, 15] in the case where the measures $\mu_n^k$, $k = 0, 1, \ldots$, are the marginals of a Markov chain on $\Omega_n$ with stationary probability $\nu_n$. See also [4] and [34], Section 2.4.

REMARKS. 1. For simplicity, we have restricted attention to the case where the maximum of the "distance" $D$ is 1. However, Definitions 2.1 and 2.2 can easily be extended to the case where $D$ is unbounded, for example, the $L^2$ distance

$$(2.1) \qquad D(\mu, \nu) = \left( \sum_{\omega \in \Omega} \left| \frac{\mu(\omega)}{\nu(\omega)} - 1 \right|^2 \nu(\omega) \right)^{1/2}.$$

In this case, in part (b) of each of the two definitions above, simply replace the limit value 1 by $\infty$ (see, e.g., [34], Section 2.4).

2. Note that the definitions above do not require that the sequence $t_n$ tends to $\infty$ (this condition is required in [15, 34], but probably for no good reason. Thanks to Jim Fill for pointing this out to us). For instance, let $\Omega_n = \{1, \ldots, n\}$, $K_n(x, y) = 1/(n-1)$ for all $y \neq x \in \Omega_n$ and $\nu_n \equiv 1/n$. This family presents a $(1, \varepsilon_n)$ total variation cut-off for any sequence $\varepsilon_n$ tending to 0. Indeed, $\|K_n^2 - \nu_n\|_{\mathrm{TV}} = 1/(n-1)$, whereas $\|K_n^0 - \nu_n\|_{\mathrm{TV}} = (n-1)/n$.

3. If a family $(\mu_n^t)$ with continuous parameter $t$ has both a $(t_n)$-cut-off and an $(s_n)$-cut-off, then $s_n \sim t_n$ (i.e., $\lim_{n \to \infty} s_n/t_n = 1$). This is also true for a discrete-time family if one of $t_n$ or $s_n$ tends to infinity. However, for a discrete-time family having both a $(t_n)$ and an $(s_n)$-cut-off, the best that can be said in general is that the limit points of the sequence $|t_n - s_n|$ all belong to the interval $[0, 1]$. Because of this, cut-off sequences that do not tend to infinity have to be treated with some special care in discrete time.

Examples of finite Markov chain cut-offs are discussed in [15] which poses the following questions. How widespread is the cut-off phenomenon for families of finite ergodic Markov chains and how can one recognize it?

It has been suggested by Yuval Peres that in some generality, these questions could be answered simply in terms of two parameters, namely, the $D$-mixing time

$$(2.2) \qquad \tau_n^D = \tau_n^D(\varepsilon) = \inf\{k : D(\mu_n^k, \nu_n) \leq \varepsilon\}$$

and an appropriately defined notion of spectral gap. Here $\varepsilon$ is a small fixed parameter (e.g., one often picks $\varepsilon = 1/4$). In the special case when the chain is reversible (hence diagonalizable with real eigenvalues in $[-1,1]$), set $\lambda_n = 1 - \beta_n$, where $\beta_n$ is the second largest eigenvalue. Peres' suggestion is that a $D$-cut-off occurs (say, in continuous time) if and only if the quantity $\lambda_n \tau_n^D$ tends to $\infty$.

It turns out that the condition $\lambda_n \tau_n^D$ tends to $\infty$ is indeed a necessary and sufficient condition for a cut-off to occur when $D$ is the $L^2$ distance defined in (2.1) (and also when $L^2$ is replaced by $L^p$ with $1 < p \leq \infty$). This is proved in [10, 11]. Earlier results in this direction are described in [34]. This means that when working with the $L^2$ distance (2.1), one can often assert that a cut-off exists without necessarily having to or being able to give a precise asymptotic of the cut-off time.

Examples due to David Aldous show that if $D$ is total variation, the condition $\lambda_n \tau_n^D \to \infty$ does not necessarily imply a cut-off (see, e.g., [10] where a version of Aldous' example is presented). Still, it is natural to conjecture that the condition that $\lambda_n \tau_n^D$ tends to infinity is sufficient under additional assumptions, for example, for random walks on finite groups based on small generating sets or for birth and death chains.

In this paper, we consider the case of continuous-time birth and death chains on $\{0, 1, \ldots, m_n\}$ started at 0. When $D$ is separation, we show that there is a cut-off if and only if $\lambda_n \tau_n^D$ tends to infinity. This will follow from previous work of Diaconis and Fill [16], who produced optimal strong stationary times through the construction of strong stationary duals, and works of Keilson [29] and Brown and Shao [8] linking spectral data to first passage times.

**3. Remarks on total variation and separation cut-offs.** For $n = 1, 2, \ldots$, let $\Omega_n$ be a finite set equipped with a Markov kernel $K_n(x, y)$ with stationary probability $\nu_n$. Fix a starting point $x_n \in \Omega_n$ and consider the sequence of probability measures $\mu_n^k$, $k = 0, 1, \ldots$, where $\mu_n^k$ is the distribution of the associated Markov chain started at $x_n$ after $k$ steps.

Let $D$ stand for either the total variation distance or separation. One thing these two notions have in common is that given the data above, there exists a sequence of real nonnegative random variables $T_n^D$ having a useful interpretation and such that

$$D(\mu_n^k, \nu_n) = P(T_n^D > k).$$

When $D$ is total variation, the $T_n^D$'s are "optimal coupling times," whereas when $D$ is separation the $T_n^D$'s are "optimal strong stationary times" (see [3, 30] and the references therein). Let $t_n$ and $\sigma_n^2$ be, respectively, the mean and variance of the random variable $T_n^D$. By a well-known form of Chebyshev's inequality (e.g., [24], (7.5), page 152), for all $a > 0$, we have

(3.1) $\quad P(T_n^D > t_n + a\sigma_n) \leq \dfrac{1}{1+a^2}, \qquad P(T_n^D < t_n - a\sigma_n) \leq \dfrac{1}{1+a^2}.$

From these facts we can draw the following conclusions:

(a) For any $\varepsilon \in (0,1)$, the mixing time $\tau_n^D(\varepsilon)$ defined in (2.2) satisfies

$$t_n - (\varepsilon^{-1} - 1)^{-1/2}\sigma_n \leq \tau_n^D(\varepsilon) \leq t_n + (\varepsilon^{-1} - 1)^{1/2}\sigma_n.$$

(b) If there is a constant $c > 0$ such that $ct_n \geq \sigma_n$ and there is cut-off at time $s_n$ with $\lim_{n\to\infty} s_n = \infty$, then $s_n \sim t_n$. When working in continuous time, the conclusion $s_n \sim t_n$ holds true without having to assume that $s_n$ tends to infinity.

(c) If $\sigma_n^{-1}t_n \to \infty$, then there is a $(t_n, \sigma_n)$-$D$-cut-off.

The upper bound in part (a) follows from the first inequality in (3.1) by solving $(1+a^2)^{-1} = \varepsilon$. The lower bound follows from the second inequality in (3.1) by solving $(1+a^2)^{-1} = 1 - \varepsilon$.

Part (b) requires a little work and we treat only the continuous-time case. Assume there is a cut-off at time $s_n$ and fix $\eta, \varepsilon \in (0,1)$. Then, for $n$ sufficiently large, we must have

$$(1-\eta)s_n \leq \tau_n^D(\varepsilon) \leq (1+\eta)s_n.$$

Setting $\varepsilon = (1+\eta^2)^{-1}$ and using the first bound in (a), we obtain

$$(1-\eta)s_n \leq t_n + \eta\sigma_n \leq (1+c\eta)t_n.$$

Using $\varepsilon = (1+\eta^{-2})^{-1}$ and the second bound in (a) gives

$$(1-c\eta)t_n \leq t_n - \eta\sigma_n \leq (1+\eta)s_n.$$

This shows that $s_n \sim t_n$, as desired.

In general, little is known about the times $T_n^D$, so these remarks have only theoretical value. In particular, we know of no nontrivial cases where an optimal coupling time has been constructed in a useful way. In contrast, there are several known examples of optimal strong stationary times to which the remarks above apply (e.g., the top to random and riffle shuffles discussed in [14] and the geometric walks in [31]).

In this context, the challenge posed by Peres' question is to relate the condition $\sigma_n^{-1}t_n \to \infty$ to spectral information. When that can be done, the remarks above may yield useful results. This will be illustrated below.

**4. Separation for birth and death chains.** Let $\Omega = \{0, \ldots, m\}$. A *birth and death chain* is a Markov chain $K$ on $\Omega$ such that $K(x,y) = 0$ unless $|x - y| \leq 1$. Write

$$\begin{aligned} q_x &= K(x, x-1), & x &= 1, \ldots, m, \\ r_x &= K(x, x), & x &= 0, \ldots, m, \\ p_x &= K(x, x+1), & x &= 0, \ldots, m-1 \end{aligned}$$

and, by convention, $q_0 = p_m = 0$. We will assume throughout that the chain is irreducible, that is, that $q_x > 0$ for $0 < x \leq m$ and $p_x > 0$ for $0 \leq x < m$. Such chains have stationary probability

$$\nu(x) = c \prod_{y=1}^{x} \frac{p_{y-1}}{q_y},$$

where $c = \nu(0)$ is a normalizing constant. Birth and death chains are in fact reversible, that is, they satisfy

$$\nu(x)K(x,y) = \nu(y)K(y,x).$$

It follows that the operator $K: L^2(\Omega, \nu) \to L^2(\Omega, \nu)$ defined by $f \mapsto Kf = \sum_y K(\cdot, y) f(y)$ is self-adjoint and thus diagonalizable with real eigenvalues in $[-1, 1]$. Let $\lambda_i$, $i = 0, \ldots, m$, be the eigenvalues of $I - K$ in nondecreasing order ($I$ denotes the identity operator). Thus, $\lambda_0 = 0 < \lambda_1 \leq \lambda_2 \leq \cdots \leq \lambda_m \leq 2$. The irreducibility of the chain is reflected in the fact that $\lambda_1 > 0$. It is also well known that $\lambda_m = 2$ if and only if the chain is periodic (of period 2) which is the case if and only if $r_x = 0$ for all $x$. In fact, because we are dealing here with irreducible birth and death chains, it is known that the $\lambda_i$'s are all distinct (e.g., [8, 29]). Karlin and McGregor [27, 28] observed that the spectral analysis of any given birth and death chain can be treated as an orthogonal polynomial problem. This sometimes leads to the exact computation of the spectrum (see, e.g., [26, 27, 28, 36] and also [32] for a somewhat different approach based on continued fractions).

Given a birth and death chain as above, let $\mu^k$ be its distribution after $k$ steps starting at 0. Let $\gamma^t$ be the distribution at time $t \geq 0$ of the associated continuous-time process started at 0, that is,

$$\gamma^t = e^{-t} \sum_0^{\infty} \frac{t^k}{k!} \mu^k.$$

In [16], Diaconis and Fill construct what they call a *strong stationary dual* for any discrete-time birth and death chain satisfying the condition $p_x + q_{x+1} \leq 1$, $0 \leq x < m$, (such chains are called *monotone chains*). The dual chain is a birth and death chain with the same eigenvalues as the original chain. The first passage time at the extremity $m$ for that dual chain is a strong stationary time for the original chain. The first passage time distribution is explicitly computed by Keilson and by Brown and Shao [8, 29] in terms of the spectral data. Fill [25] treats continuous-time chains (the condition $p_x + q_{x+1} \leq 1$, $0 \leq x < m$, is not needed in that case). These works give the following result:

THEOREM 4.1 ([16, 25]). *Let $K$ be an irreducible birth and death chain as above.*

(a) *For the associated continuous-time process started at 0, we have*

$$\mathrm{sep}(\gamma^t, \nu) = \max_{0 \leq x \leq m} \left\{ 1 - \frac{\gamma^t(x)}{\nu(x)} \right\} = \sum_{i=1}^m \prod_{j \neq i} \frac{\lambda_j}{\lambda_i - \lambda_j} e^{-t\lambda_i}.$$

(b) *For the discrete-time chain, assuming that $p_x + q_{x+1} \leq 1$, $0 \leq x < m$, we have*

$$\mathrm{sep}(\mu^k, \nu) = \max_{0 \leq x \leq m} \left\{ 1 - \frac{\mu^k(x)}{\nu(x)} \right\} = \sum_{i=1}^m \prod_{j \neq i} \frac{\lambda_j}{\lambda_i - \lambda_j} (1 - \lambda_i)^k.$$

Although these are beautiful formulas, it is not so obvious how to use them to derive explicit bounds. However, (a) has a very clear interpretation: it says that separation at time $t$ is the tail of a sum of $m$ independent exponential random variables with respective parameters $\lambda_i$, $1 \leq i \leq m$. Similarly, when all $\lambda_i$ are in $[0, 1]$, (b) says that separation at time $k$ is the tail of a sum of independent geometric random variables with respective parameters $\lambda_i$, $1 \leq i \leq m$. In particular, we have the following obvious corollary:

COROLLARY 4.2. *Let $K$ be an irreducible birth and death chain as above.*

(a) *For the associated continuous-time process started at 0, we have*

$$\mathrm{sep}(\gamma^t, \nu) = \max_{0 \leq x \leq m} \left\{ 1 - \frac{\gamma^t(x)}{\nu(x)} \right\} = P(T > t)$$

*where $T = \sum_1^m S_i$, each $S_i$ being an exponential random variable with parameter $\lambda_i$ and the $S_i$'s being independent. In particular,*

$$E(T) = \sum_1^m \lambda_i^{-1}, \qquad \mathrm{Var}(T) = \sum_1^m \lambda_i^{-2}.$$

(b) *For the discrete-time chain, assuming that $p_x + q_{x+1} \leq 1$, $0 \leq x < m$, we have*

$$\mathrm{sep}(\mu^k, \nu) = \max_{0 \leq x \leq m} \left\{ 1 - \frac{\mu^k(x)}{\nu(x)} \right\} = P(T > k),$$

*where $T$ is a random variable with*

$$E(T) = \sum_1^m \lambda_i^{-1}, \qquad \mathrm{Var}(T) = \sum_1^m (1 - \lambda_i) \lambda_i^{-2}.$$

*The random variable $T$ can be written as a sum $T = \sum_1^m S_i$ where the random variables $S_i$, $1 \leq i \leq m$, are independent and $S_i$ is geometric with probability of success $\lambda_i$ if $\lambda_i \in (0, 1]$, whereas $S_i$ is a Bernoulli variable with parameter $\lambda_i^{-1}$ if $\lambda_i > 1$.*

REMARKS. 1. The times $S_i$ have no known interpretations in terms of the underlying birth and death chain.

2. Of course, the same results apply if the birth and death chain starts at the other extremity $m$. As the spectral data does not change, it follows that $\operatorname{sep}(\gamma^t(0,\cdot),\nu)$, that is, the separation starting from 0, and $\operatorname{sep}(\gamma^t(m,\cdot),\nu)$, that is, the separation starting from $m$, are equal at all times! This is in sharp contrast with what happens in total variation distance, for which starting at one or the other extremity can lead to very different behaviors.

3. In view of the above results, and from the viewpoint developed in the next few sections, it is interesting to note that for any set of $m$ distinct positive numbers $0 < \lambda_1 < \cdots < \lambda_m \leq 1$, there is a birth and death chain as above with eigenvalues $(0,\lambda_1,\ldots,\lambda_m)$; see [33] and the references therein.

4. Recall that, given an arbitrary positive distribution $\nu$ on $\{0,1,\ldots,m\}$, the Metropolis algorithm based on symmetric random walk on the path $\{0,\ldots,m\}$ produces an irreducible birth and death chain having $\nu$ as its stationary distribution. Thus, any positive distribution can occur as the stationary distribution of an irreducible birth and death chain.

5. Continuous-time Markov chains offer the freedom to choose a time scale. Starting with a Markov kernel $K(x,y)$ on a countable space, we can consider the continuous-time Markov chain generated by the $K - I$, where $I$ denotes the identity matrix. Starting at $x$, this continuous-time chain has probability distribution at time $t$ given by

$$\gamma^t(x,\cdot) = e^{-t} \sum \frac{t^n}{n!} K^n(x,\cdot).$$

However, the most general and natural definition of a countable continuous-time Markov chain involves a matrix $Q(x,y)$ (representing the generator) satisfying $\sum_y Q(x,y) = 0$ and $Q(x,y) \geq 0$ if $x \neq y$. In that generality, the quantity $\sum_{y \neq x} Q(x,y)$ does not have to be uniformly bounded (and explosion in finite time is possible). On a finite state space, we can always set $q = \max_x\{-Q(x,x)\}$ and consider the (discrete-time) chain with kernel $K(x,y) = I(x,y) + q^{-1}Q(x,y)$ where $I$ is the identity matrix. Letting $\gamma_Q^t(x,\cdot)$ be the probability distribution of the continuous-time Markov chain with generator $Q$ started at $x$, we have

$$\gamma_Q^t(x,\cdot) = \gamma^{qt}(x,\cdot).$$

Let us now consider a family of continuous-time (finite state space) ergodic Markov chains $(\Omega_n, \gamma_{Q_n}^t(x_n,\cdot), \nu_n)$ and consider whether or not this family presents a $D$-cut-off. The answer to this question is independent of the chosen time scale. Indeed, using the notation introduced above,

$$\gamma_{Q_n}^t(x_n,\cdot) = \gamma_n^{q_n t}(x_n,\cdot).$$

It follows that $(\Omega_n, \gamma_n^t(x_n, \cdot), \nu_n)$ presents a $D$-cut-off at time $t_n$ [resp. a $D$-cut-off of type $(t_n, b_n)$] if and only if $(\Omega_n, \gamma_{Q_n}^t(x_n, \cdot), \nu_n)$ presents a $D$-cut-off at time $t_n/q_n$ [resp. a $D$-cut-off of type $(t_n/q_n, b_n/q_n)$]. This remark would not be valid if we had required that $(t_n)$ must tend to infinity in the definition of a cut-off.

**5. Separation cut-off for birth and death chains.** We now describe what the previous section entails concerning the cut-off phenomenon. For $n = 1, 2, \ldots$, let $\Omega_n = \{0, 1, \ldots, m_n\}$ be equipped with an irreducible birth and death chain $K_n$ having stationary measure $\nu_n$. Let $q_{n,x}, r_{n,x}$ and $p_{n,x}$ be the corresponding transition probabilities.

Let $\mu_n^k$ be the distribution of the associated chain at time $k$ started at 0. Let $\gamma_n^t$ be the distribution of the continuous-time process at time $t$ started at 0. Let $\lambda_{n,i} \in [0, 2]$, $0 \le i \le m_n$, be the corresponding eigenvalues. Set

$$\lambda_n = \lambda_{n,1}, \qquad t_n = \sum_1^{m_n} \lambda_{n,i}^{-1}.$$

Finally, for any $\varepsilon \in (0, 1)$, consider the separation mixing time

$$\tau_n(\varepsilon) = \inf\{t : \operatorname{sep}(\gamma_n^t, \nu_n) \le \varepsilon\}.$$

THEOREM 5.1. *Referring to the setting and notation introduced above, the family*

$$(\Omega_n, \nu_n, (\gamma_n^t)_{t>0})_{n=1,2,\ldots}$$

*has a separation cut-off if and only if $N_n = \lambda_n t_n$ tends to infinity. For any $c > 0$, the separation bounds*

$$(5.1) \quad \operatorname{sep}(\gamma_n^{(1+c)t_n}, \nu_n) \le \frac{1}{1 + c^2 N_n}, \qquad \operatorname{sep}(\gamma_n^{(1-c)t_n}, \nu_n) \ge 1 - \frac{1}{1 + c^2 N_n}$$

*always hold and for any fixed $\varepsilon \in (0, 1)$, the condition $\lambda_n t_n \to \infty$ is equivalent to $\lambda_n \tau_n(\varepsilon) \to \infty$.*

PROOF. By Corollary 4.2, we have

$$\operatorname{sep}(\mu_n^t, \nu_n) = P(T_n > t),$$

where $T_n$ has mean $t_n$ and variance $\operatorname{Var}(T_n) = \sigma_n^2$ satisfying

$$\sigma_n^2 = \sum_1^{m_n} \lambda_{n,i}^{-2} = \lambda_n^{-2} \sum_1^{m_n} (\lambda_n/\lambda_{n,i})^2$$

$$\le \lambda_n^{-2} \left( \sum_1^{m_n} \lambda_n/\lambda_{n,i} \right) = \lambda_n^{-1} t_n.$$

Here we have simply used the fact that $\lambda_n/\lambda_{n,i} \leq 1$ to obtain the middle inequality. Hence, we have

(5.2) $$\sigma_n \leq t_n \quad \text{and, better,} \quad \sigma_n \leq N_n^{-1/2} t_n.$$

The separation bounds (5.1) follow directly from (5.2) and the Chebyshev inequalities (3.1).

Assume that $N_n = \lambda_n t_n \to \infty$. By the second inequality in (5.2), it follows that $t_n/\sigma_n \to \infty$. By (c) of Section 3, there is a separation cut-off at time $t_n$ and, even better, a $(t_n, \sigma_n)$-cut-off. Conversely, if there is a cut-off at time $s_n$, then by (5.2) and (b) of Section 3, we must have $s_n \sim t_n$ and there must be a cut-off at time $t_n$. By (5.1), this implies that $N_n$ tends to infinity.

Now, fix $\varepsilon \in (0,1)$. By the upper bound in (a) of Section 3 and the first inequality in (5.2), we have

(5.3) $$\tau_n(\varepsilon) \leq t_n + (\varepsilon^{-1}-1)^{1/2}\sigma_n \leq (1+(\varepsilon^{-1}-1)^{1/2})t_n.$$

Hence $\lambda_n \tau_n(\varepsilon) \to \infty$ implies $\lambda_n t_n \to \infty$. Conversely, if $t_n \lambda_n \to \infty$, then there is a cut-off at time $t_n$ and by (5.2) and (a) of Section 3, $t_n \sim \tau_n(\varepsilon)$. It follows that $\lambda_n \tau_n(\varepsilon) \to \infty$. This ends the proof of Theorem 5.1. $\square$

REMARKS. 1. Theorem 5.1 shows that for continuous-time birth and death chains started at 0, a separation cut-off can occur only if $m_n$ tends to infinity.

2. For $D$ and $(\Omega_n, \nu_n, (\mu_n^k)_{k=0,1,...})$ as in Definition 2.1, we say that there is a *D-precut-off* at time $s_n$ if there are constants $0 < c \leq 1 \leq C < \infty$ such that

$$\lim_{n\to\infty} D(\mu_n^{k_n}, \nu_n) \to \begin{cases} 0, & \text{if } k_n \geq Cs_n, \\ 1, & \text{if } k_n \leq cs_n. \end{cases}$$

Obvious modifications apply in continuous time. Theorem 5.1 shows that there cannot be a separation precut-off if $\lambda_n t_n$ is bounded. Hence, for continuous birth and death chains started at 0, the existence of a separation precut-off is equivalent to the existence of a separation cut-off.

3. In (5.1), we can replace $N_n = \lambda_n t_n$ by $N'_n = \lambda_n \tau_n(\varepsilon)$ as follows. The bound (5.3) and the first inequality in (5.2) give $N'_n \leq (1+(\varepsilon^{-1}-1)^{1/2})N_n$. Hence, for $c > 0$ and $\varepsilon \in (0,1)$, we have

$$\text{sep}(\gamma^{(1+c)t_n}, \nu_n) \leq \frac{1}{1+c^2(1+(\varepsilon^{-1}-1)^{1/2})^{-1}N'_n}$$

and

$$\text{sep}(\gamma^{(1-c)t_n}, \nu_n) \geq 1 - \frac{1}{1+c^2(1+(\varepsilon^{-1}-1)^{1/2})^{-1}N'_n}.$$

The same remark applies to (6.3) below which improves upon (5.1).

The next result is the discrete-time version of Theorem 5.1. It requires the "monotonicity" assumption $p_x + q_{x+1} \leq 1$. The proof is similar to that of Theorem 4.2(b) and is therefore omitted.

THEOREM 5.2. *Referring to the setting and notation introduced above, assume that for each $n$ and each $x \in \{0, \ldots, m_n - 1\}$, we have*

$$p_{n,x} + q_{n,x+1} \leq 1.$$

*Then the family*

$$(\Omega_n, \nu_n, (\mu_n^k)_{k=0,1,\ldots})_{n=1,2,\ldots}$$

*has a separation cut-off if and only if $N_n = \lambda_n t_n$ tends to infinity.*

The so-called monotonicity condition $p_{n,x} + q_{n,x+1} < 1$ easily implies that $r_{n,0} > 0$, thus ensuring aperiodicity. It is, however, a little surprising that negative eigenvalues of $K$ (i.e., $1 - \lambda_i$ with $\lambda_i > 1$) play no role whatsoever in Theorem 5.2. As in the continuous case, for any fixed $\varepsilon \in (0,1)$, the theorem above can be stated using $\tau_n(\varepsilon) = \inf\{k : \text{sep}(\mu_n^k, \nu_n) \leq \varepsilon\}$ instead of $t_n$.

**6. The shape of the cut-off.** When a cut-off is determined, say at time $s_n$, the next task is to look at the window size. If one is able to establish an $(s_n, b_n)$-cut-off [possibly adjusting the sequence $(s_n)$], then the question of the optimality of the window size $b_n$ is posed. One way to answer this question is to obtain the shape of the cut-off, that is, to determine the functions $f_\pm$ of Definition 2.2. If $f_+$ and $1 - f_-$ are nonzero in a neighborhood of $\infty$, then the sequence $(b_n)$ is optimal and the functions $f_\pm$ describe the shape of the cut-off. Only a small number of such results have been established (see, e.g., [15]). In the cases of interest to us in this paper, Corollary 4.2 easily allows us to obtain the shape of the cut-off.

THEOREM 6.1. *Referring to a family of birth and death chains as in Section 5 and using the notation introduced there, assume that $N_n = \lambda_n t_n \to \infty$ and set*

$$\sigma_n^2 = \sum_1^{m_n} \lambda_{n,i}^{-2}.$$

(a) *Assume that $\lambda_n \sigma_n \to \infty$. Then for any real $c$,*

$$\lim_{n \to \infty} \text{sep}(\gamma_n^{t_n + c\sigma_n}, \nu_n) = 1 - \Phi(c) \qquad \text{where } \Phi(t) = \frac{1}{\sqrt{2\pi}} \int_{-\infty}^t e^{-x^2/2} \, dx.$$

*In particular, there is a $(t_n, \sigma_n)$-cut-off, but no $(t_n, \lambda_n^{-1})$-cut-off.*

(b) *Assume that $\lambda_n \sigma_n$ is bounded. Then there is a $(t_n, \sigma_n)$-cut-off [equivalently a $(t_n, \lambda_n^{-1})$-cut-off] and for any real $c > 0$, we have*

$$\liminf_{n \to \infty} \operatorname{sep}(\gamma_n^{t_n + c\sigma_n}, \nu_n) > 0,$$

*whereas for any real $c < 0$,*

$$\limsup_{n \to \infty} \operatorname{sep}(\gamma_n^{t_n + c\sigma_n}, \nu_n) < 1.$$

PROOF. We have $\operatorname{sep}(\gamma_n^t, \nu_n) = P(T_n > t)$. Consider the moment generating function

$$M_n(t) = E(e^{t(T_n - t_n)/\sigma_n}).$$

As $T_n$ is a sum of $m_n$ independent exponential random variables with respective parameters $\lambda_{n,i}$, $i = 1, \ldots, m_n$, we have

$$M_n(t) = e^{-tt_n/\sigma_n} \prod_{i=1}^{m_n} \frac{\lambda_{n,i}}{\lambda_{n,i} - t/\sigma_n} = e^{F_n(t)},$$

where

$$F_n(t) = -t_n \sigma_n^{-1} t - \sum_{i=1}^{m_n} \log(1 - \lambda_{n,i}^{-1} \sigma_n^{-1} t) = \sum_{k=2}^{\infty} \frac{\theta_k(n)}{k \theta_2(n)^{k/2}} t^k$$

with

$$\theta_k(n) = \sum_{i=1}^{m_n} (\lambda_{n,1}/\lambda_{n,i})^k.$$

As $\lambda_{n,1}/\lambda_{n,i} \leq 1$, we have $1 \leq \theta_k(n) \leq \theta_2(n), k \geq 2$. Hence, the series above converges, at least for $t \in (-1, 1)$ and

(6.1) $$0 \leq F_n(t) - t^2/2 \leq \sum_{k=3}^{\infty} \frac{t^k}{k \theta_2(n)^{(k-2)/2}}.$$

If $\lambda_n \sigma_n = \theta_2(n)^{1/2} \to \infty$, then $M_n(t)$ tends to $e^{t^2/2}$ for any real $t$ and it follows that $\sigma_n^{-1}(T_n - t_n)$ is asymptotically distributed as a standard normal random variable; see [12]. This proves part (a) of Theorem 6.1.

Assume now that $\lambda_n \sigma_n$ is bounded, say $\lambda_n \sigma_n \leq A$. Hence,

$$\lambda_n^{-1} \leq \sigma_n \leq A \lambda_n^{-1}$$

and, for any $k = 2, 3, \ldots,$

$$1 \leq \theta_k(n) \leq A.$$

Obviously (e.g., by Chebyshev's inequality), the distributions of $\sigma_n^{-1}(T_n - t_n)$, $n = 1, 2, \ldots,$ form a tight family. Given any subsequence $(n_j)$, we can extract from it a subsequence $(n_{j_\ell})$ such that along that subsequence, $P(T_n >$ $t_n + c\sigma_n$) converges to $P(T > c)$, $c \in \mathbb{R}$, for some random variable $T$. Now, from the previous moment generating function computation, it follows that along $(n_{j_\ell})$, the limit

$$\lim_{\ell \to \infty} \theta_k(n_{j_\ell}) = \theta_k \in [1, A]$$

exists for each $k \geq 2$ and $T$ has moment generating function

(6.2) $$\exp\left(\frac{t^2}{2} + \sum_{k \geq 3} \frac{\theta_k t^k}{k \theta_2^{k/2}}\right)$$

for all $t \in (-\theta_2^{1/2}, \theta_2^{1/2})$. As the variables $\sigma_n^{-1}(T_n - t_n)$ are infinitely divisible, $T$ is also infinitely divisible (see, e.g., [35] for properties of infinitely divisible distributions). By (6.2), the normal component of $T$ is nontrivial and it follows that $T$ has a smooth positive density. Obviously, this implies part (b) of Theorem 6.1. Note that no limit points of the sequence $(T_n - t_n)/\sigma_n$ can be normal. This ends the proof of Theorem 6.1. □

Let us observe that the first part of Theorem 6.1(b), that is, the fact that for any real $c > 0$,

$$\liminf_{n \to \infty} \operatorname{sep}(\gamma_n^{t_n + c\sigma_n}, \nu_n) > 0,$$

can be proved by a very elementary argument. To bound $P(T_n > t_n + c\sigma_n)$ from below, write

$$P(T_n > t_n + c\sigma_n)$$
$$\geq P\left(S_{n,1} > \lambda_n^{-1} + (c+1)\sigma_n; \sum_{i=2}^{m_n} S_{n,i} > \sum_{i=2}^{m_n} \lambda_{n,i}^{-1} - \sigma_n\right)$$
$$\geq P(S_{n,1} > \lambda_n^{-1} + (c+1)\sigma_n) P\left(\sum_{i=2}^{m_n} S_{n,i} > \sum_{i=2}^{m_n} \lambda_{n,i}^{-1} - \sigma_n\right).$$

By the second inequality in (3.1), and the fact that $\sigma_n^2 \geq \operatorname{Var}(\sum_{i=2}^{m_n} S_{n,i})$, we have

$$P\left(\sum_{i=2}^{m_n} S_{n,i} > \sum_{i=2}^{m_n} \lambda_{n,i}^{-1} - \sigma_n\right) \geq \tfrac{1}{2}.$$

It follows that

$$P(T_n > t_n + c\sigma_n) \geq \tfrac{1}{2} e^{-(\lambda_n^{-1} + (c+1)\sigma_n)\lambda_n} \geq \tfrac{1}{2} e^{-(1+(c+1)A)}.$$

In contrast, the inequality concerning the lower tail seems harder to prove without the sophisticated tools of infinitely divisible distributions.

REMARKS. 1. In part (b) of Theorem 6.1, assume further that for each $k$,

$$\theta_k = \lim_{n\to\infty} \left(\frac{\lambda_n}{\lambda_{n,i}}\right)^k < \infty$$

exists. Then for any $c \in \mathbb{R}$, we have

$$\lim_{n\to\infty} \text{sep}(\gamma_n^{t_n+c\sigma_n}, \nu_n) = 1 - F(c),$$

where $F(t)$ is the distribution function of an infinitely divisible law whose moment generating function is given in $(-\theta_2^{-1/2}, \theta_2^{-1/2})$ by (6.2). In particular, $0 < F(c) < 1$ for all $c \in \mathbb{R}$. For instance, the Bernoulli–Laplace example in Section 7 has a nonnormal cut-off shape in separation.

2. It follows from (6.1) that $E(e^{\pm(T_n-t_n)/(2\sigma_n)}) \leq e^{1/4}$. Hence, for all $c > 0$,

$$P(T_n > t_n + c\sigma_n) \leq e^{-(c-1/2)/2}, \qquad P(T_n > t_n - c\sigma_n) \geq 1 - e^{-(c-1/2)/2}$$

and thus

(6.3)
$$\text{sep}(\gamma^{(1+c)t_n}, \nu_n) \leq e^{-(cN_n-1/2)/2},$$
$$\text{sep}(\gamma^{(1-c)t_n}, \nu_n) \geq 1 - e^{-(cN_n-1/2)/2},$$

where $N_n = \lambda_n t_n$ tends to infinity with $n$.

**7. Examples.** This section illustrates our results by looking at various explicit (and not so explicit) families of birth and death chains.

*Simple random walk.* Consider the simple random walk on $\{0, \ldots, n\}$ with $r_0 = r_n = p_0 = p_j = q_j = q_n = 1/2$, $j \in \{1, \ldots, n-1\}$. We refer the reader to [24], XVI.3, for the spectral information used below. We have

$$\lambda_{n,j} = 1 - \cos\frac{\pi j}{n+1}, \qquad j = 0, \ldots, n.$$

As $\lambda_{n,j} \geq (j/(n+1))^2$ (and this is optimal, up to a multiplicative constant), we see that $\lambda_n^{-1} = \lambda_{n,1}^{-1}$, $t_n = \sum \lambda_{n,j}^{-1}$ and $\sigma_n = (\sum \lambda_{n,j}^{-2})^{1/2}$ are all of order $n^2$. By Theorems 5.1 and 5.2, there is no separation cut-off (either in continuous or discrete time). Of course, this is well known!

Suppose, instead, that $r_n = p_j = p$, $0 \leq j \leq n-1$ and $q_j = r_0 = q$, $1 \leq j \leq n$, $p+q = 1$, $0 \leq q \leq p \leq 1$. Then the eigenvalues are

$$\lambda_{n,j} = 1 - 2\sqrt{pq}\cos\frac{\pi j}{n+1}, \qquad j = 0, \ldots, n.$$

Hence, if $p$ and $q$ are fixed with $0 < q < p < 1$, we have $\lambda_n = 1 - 2\sqrt{pq}$, whereas $t_n = \sum_1^n \lambda_{n,j}^{-1}$ is greater than $n/2$ and $\sigma_n$ is greater than $\sqrt{n/2}$.

Theorems 5.1, 5.2 and 6.1 prove the existence of a $(t_n, \sigma_n)$ separation cut-off with a normal shape. Observe that

$$t_n = an + O(1), \qquad a = \int_0^1 (1 - 2\sqrt{pq}\cos(\pi x))^{-1}\, dx = \frac{1}{\sqrt{1-4pq}}$$

and

$$\sigma_n^2 = bn + O(1), \qquad b = \int_0^1 (1 - 2\sqrt{pq}\cos(\pi x))^{-2}\, dx.$$

Note that the window of the separation cut-off is not given by $\lambda_n^{-1} \simeq 1$ in this case.

Using diagonalization, one finds that (starting at 0) this chain has a $(2cn, 1)$-$L^\infty$ cut-off and a $(cn, 1)$-$L^2$ cut-off with

$$c = \frac{\log(p/q)}{2(1 - 2\sqrt{pq})}.$$

It is a calculus exercise to check that $c > a$, that is,

$$\log(p/q) \geq \frac{2(1 - 2\sqrt{pq})}{\sqrt{1-4pq}}.$$

Indeed, writing $p = (1+u)/2$, $q = (1-u)/2$, $u \in (0,1)$, we get $1 - 4pq = u^2$ and the above inequality boils down to

$$\frac{\log(1+u) - \log(1-u)}{2u} \geq \frac{1 - \sqrt{1-u^2}}{u^2},$$

which holds true for $u \in (0,1)$ because the left-hand side is at least 1, whereas the right-hand side is at most 1. Thus, the $L^2$-cut-off occurs later than the separation cut-off (this is not always true—there are many examples where the separation cut-off is twice the $L^2$-cut-off).

Diaconis and Fill [16], Example 4.46, and Belsley [6], Chapter V, study various versions of this chain in detail (in discrete time) and show that there is a cut-off in total variation at time $an$ with an optimal window of size $\sqrt{n}$. The fact that the total variation cut-off is the same as the separation cut-off can be explained as follows. If the chain starts from the top point $n$, it is not hard to use the available spectral information to show that it converges in a constant number of steps (most of the mass is around $n$). To understand total variation starting from the bottom point 0, it will thus be enough to analyze the first time one hits $n$. By [8, 29], this first hitting time is equal in law to the optimal strong stationary time of Corollary 4.2(a).

Note the very different window sizes, namely, of order 1 for the $L^2$ and $L^\infty$ cut-offs and of order $\sqrt{n}$ for the separation and total variation cut-off.

*Metropolis chains.* Now consider an arbitrary probability distribution $\nu$ on $\{0,\ldots,n\}$ with $\nu(j) > 0$, $j \in \{0,\ldots,n\}$. Use the Metropolis algorithm with base chain the simple symmetric random walk above to obtain a birth and death chain with stationary measure $\nu$ (see, e.g., [20]). By construction, this chain satisfies the monotonicity condition $p_x + q_{x+1} \leq 1$ and Theorems 5.1, 5.2 and 6.1 apply.

For instance, if $\nu(j) = a(1+j)^d$, then the results in [20] show that there is no cut-off (in total variation or separation) and that $\lambda_n$ is of order $n^{-2}$. Hence, it follows from Theorem 5.1 that $t_n = \sum_1^n \lambda_{n,i}^{-1}$ must be of order $n^2$, that is, the eigenvalues $\lambda_i$ must grow rapidly enough from their minimum of order $n^{-2}$. We do not know if this can be easily checked by bounding higher eigenvalues. For instance, it does not follow from the (rather sophisticated) eigenvalue bound $\lambda_{n,i} \geq c i^{2/d} n^{-2}$ given by [34], Theorem 3.4.4, and [20].

As a second example, take $\nu(j) = 2^{-n}\binom{n}{j}$. It is proved in [21] that, for this example, $\lambda_n$ is of order $1/n$ and it follows from the proof that in fact $\lambda_{n,i}$ is of order $i/n$. Hence, $t_n$ is of order $n \log n$, $\lambda_n t_n \to \infty$ and $\lambda_n \sigma_n$ is bounded. By Theorems 5.1 and 5.2, there is a $(t_n, \lambda_n^{-1})$ cut-off in separation (starting from 0). The exact asymptotic behavior of the cut-off time $t_n$ is not known.

*Bernoulli–Laplace models.* Consider two urns, the left containing $r$ red balls and the right containing $(n-r)$ black balls, with $0 < 2r \leq n$. At each step, a ball is picked uniformly at random in each urn and the two balls are switched. The process is completely determined by the number of red balls in the right urn and this is a birth and death chain on $\{0,\ldots,r\}$. The stationary distribution is

$$\nu_{n,r}(j) = \frac{\binom{r}{j}\binom{n-r}{r-j}}{\binom{n}{r}}$$

and for $x \in \{0,\ldots,r\}$, the rates are given by

$$p_x = \frac{(r-x)(n-r-x)}{r(n-r)}, \qquad q_x = \frac{x^2}{r(n-r)}, \qquad r_x = 1 - p_x - q_x.$$

The eigenvalues of this chain are well known (this goes back at least to [27]; see, e.g., [23]) and given by (with an obvious change in notation)

$$\lambda_{n,r,i} = \frac{i(n-i+1)}{r(n-r)}.$$

Hence, the smallest nonzero eigenvalue is

$$\lambda_{n,r} = \frac{n}{r(n-r)}$$

and we have

$$t_{n,r} = r(n-r) \sum_1^r \frac{1}{i(n-i+1)}$$
$$= \frac{r(n-r)}{n} \sum_1^r \left(\frac{1}{i} + \frac{1-1/i}{n(1-(i-1)/n)}\right)$$
$$= \frac{r(n-r)}{n}(\log r + O(1))$$

and

$$\theta_2(n,r) = \lambda_{n,r}^{-2}\sigma_{n,r}^2 = \sum_1^r \left(\frac{n}{i(n-i+1)}\right)^2 = O(1).$$

In both cases, the $O(1)$ is uniform for all $r \leq n/2$ as $r$ tends to infinity. Given this data, Theorem 5.1 shows that for any sequence $(r_\ell, n_\ell)$ with $r_\ell \to \infty$ and $r_\ell < n_\ell/2$, the associated continuous-time chain has a $(s_\ell, \xi_\ell)$-separation cut-off with

$$s_\ell = (1 - r_\ell/n_\ell)r_\ell \log r_\ell, \qquad \xi_\ell = (1 - r_\ell/n_\ell)r_\ell.$$

If $r_\ell/n_\ell$ tends to zero, then

$$t_{n_\ell, r_\ell} = r_\ell(\log r_\ell + \gamma + o(1)) \quad \text{and} \quad \theta_2(n_\ell, r_\ell) = \frac{\pi^2}{6} + o(1),$$

where $\gamma$ denotes the well-known $\gamma$ constant. In this case, a slight variation on Theorem 6.1 shows that the limit shape for the $(s_\ell, \xi_\ell)$-separation cut-off is given by the Gumbel distribution [density $\exp(-(x + e^{-x}))$ on $\mathbb{R}$].

The results above should be compared with those of [23], Theorem 2, where an $L^2$-cut-off of type $(\zeta_\ell, \xi_\ell)$ is proved with $\xi_\ell$ as above and

$$\zeta_\ell = \frac{1}{2}(1 - r_\ell/n_\ell)r_\ell \log n_\ell = \left(\frac{\log n_\ell}{2 \log r_\ell}\right)s_\ell.$$

We now describe what happens in total variation for this family. This is briefly discussed in [7], Section 1.5. Recall that the total variation distance is easily bounded by both separation and $L^2$ (see, e.g., [3, 4]). Hence, if there is a total variation cut-off, it is bounded above by

$$\rho_\ell = \min\{s_\ell, \zeta_\ell\} = \begin{cases} s_\ell, & \text{if } n_\ell \geq r_\ell^2, \\ \zeta_\ell, & \text{if } n_\ell \leq r_\ell^2 \leq n_\ell^2/4, \end{cases}$$

that is, the minimum of the separation cut-off time $s_\ell$ and the $L^2$ cut-off time $\zeta_\ell$. It turns out that this upper bound is sharp and that there is, in fact, a $(\rho_\ell, \xi_\ell)$-cut-off in total variation. We find this phenomenon quite interesting and surprising. For this natural family of chains, separation and $L^2$ cut-off times cross each other [as functions of the parameters $(n,r)$] and the total variation cut-off time is given by the minimum.

To prove this, we only need a lower bound on total variation matching the upper bound provided by the separation and $L^2$ results. This lower bound can be obtained by the method introduced in [23] and used there in the case $r = n/2$. Namely, to lower bound total variation between $\gamma_{n,r}^t$ and its stationary measure $\nu_{n,r}$, use a set of the form $A = \{\phi_1 \leq \alpha\}$, where $\phi_1$ is an eigenfunction associated with the lowest nonzero eigenvalue $\lambda_{n,r} = \lambda_{n,r,1}$. A complete set of eigenfunctions $\{\phi_i : i = 0, \ldots, r\}$ ($\phi_i = \phi_{n,r,i}$ associated with $\lambda_{n,r,i}$) is described in [23]. In particular, we can take

$$\phi_0(x) = 1, \qquad \phi_1(x) = 1 - \frac{xn}{n(n-r)}$$

and

$$\phi_2(x) = 1 - \frac{2x(n-1)}{n(n-r)} + \frac{(n-1)(n-2)x(x-1)}{r(n-r)(n-r-1)(r-1)},$$

where $x \in \{0, \ldots, r\}$. Given this data, one checks that

$$\phi_1^2 = \frac{1}{n-1}\phi_0 + \frac{n^2 - 4r(n-r)}{r(n-r)(n-2)}\phi_1 + \frac{n^2(n-r-1)(r-1)}{r(n-r)(n-1)(n-2)}\phi_2.$$

This formula allows us to compute the variance of $\phi_1$ under $\gamma_{n,r}^t$ [the variance of $\phi_1$ under the stationary measure $\nu_{n,r}$ is $(n-1)^{-1}$]. Namely,

$$\begin{aligned}
\operatorname{Var}_{\gamma_{n,r}^t}(\phi_1) &= \frac{1}{n-1} + \frac{n^2 - 4r(n-r)}{r(n-r)(n-2)} e^{-t\lambda_{n,r,1}} \\
&\quad + \frac{n^2(n-r-1)(r-1)}{r(n-r)(n-1)(n-2)} e^{-t\lambda_{n,r,2}} - e^{-2t\lambda_{n,r,1}} \\
&= \frac{1}{n-1} + \frac{n^2 - 4r(n-r)}{r(n-r)(n-2)} e^{-t\lambda_{n,r,1}} \\
&\quad + \left(\frac{n^2(n-r-1)(r-1)}{r(n-r)(n-1)(n-2)} - 1\right) e^{-t\lambda_{n,r,2}} \\
&\quad + (1 - e^{-t(2\lambda_{n,r,1} - \lambda_{n,r,2})}) e^{-t\lambda_{n,r,2}} \\
&= \frac{1}{n-1} + \frac{n^2 - 4r(n-r)}{r(n-r)(n-2)} e^{-t\lambda_{n,r,1}} \\
&\quad + \frac{r(n-r)(2-3n) + n^2(n-1)}{r(n-r)(n-1)(n-2)} e^{-t\lambda_{n,r,2}} \\
&\quad + (1 - e^{-t(2\lambda_{n,r,1} - \lambda_{n,r,2})}) e^{-t\lambda_{n,r,2}} \\
&\leq C\left(\frac{1}{n-1} + \frac{1}{r} e^{-t\lambda_{n,r,1}} + \frac{t}{r(n-r)} e^{-t\lambda_{n,r,2}}\right).
\end{aligned}$$

For the last term, we have used $1 - e^{-u} \leq u$, $u \geq 0$, and

$$2\lambda_{n,r,1} - \lambda_{n,r,2} = \frac{1}{r(n-r)}.$$

That the difference between $2\lambda_{n,r,1}$ and $\lambda_{n,r,2}$ is small is what makes this proof work.

Now we consider the two cases $r^2 \leq n$ and $n \leq r^2 \leq n^2/4$. In the first case ($r^2 \leq n$), set

$$t = \frac{1}{\lambda_{n,r,1}}(\log r - c), \qquad 0 < c < \log r,$$

and consider $A = \{\phi_1 \leq \alpha\}$ with $2\alpha = e^{-t\lambda_{n,r,1}}$. Chebyshev's inequality gives

$$\nu_{n,r}(A) - \gamma_{n,r}^t(A) \geq 1 - \frac{r^2 e^{-2c}}{n-1} - C'\left(\frac{1}{n-1} + \frac{1}{r^2}e^c + \frac{\log r}{r^2 n}e^{2c}\right)r^2 e^{-2c}$$

$$\geq 1 - C''\left(e^{-c} + \frac{\log n}{n}\right).$$

This, together with the earlier separation result, proves the existence of a $(s_\ell, \xi_\ell)$ total variation cut-off when $r_\ell^2 \leq n_\ell^2$ and $r_\ell$ tends to infinity.

In the second case ($n \leq r^2 \leq n^2/4$), set

$$t = \frac{1}{2\lambda_{n,r,1}}(\log n - c), \qquad 0 < c < \log n,$$

and consider $A = \{\phi_1 \leq \alpha\}$ with, again, $2\alpha = e^{-t\lambda_{n,r,1}}$. Now, Chebyshev's inequality gives

$$\nu_{n,r}(A) - \gamma_{n,r}^t(A) \geq 1 - \frac{n e^{-2c}}{n-1} - C'\left(\frac{1}{n-1} + \frac{1}{r\sqrt{n}}e^c + \frac{\log n}{n^2}e^{2c}\right)n e^{-2c}$$

$$\geq 1 - C''\left(e^{-c} + \frac{\log n}{n}\right).$$

This, together with the $L^2$ result, proves the existence of a $(\zeta_\ell, \xi_\ell)$ total variation cut-off when $n_\ell \leq r_\ell^2 \leq n_\ell^2/4$ and $r_\ell$ tends to infinity.

*The hypercube and Hamming chains.* Consider the set $\{0, \ldots, n-1\}^r$ and the Markov chain that picks one of the $r$ coordinates uniformly at random and changes this coordinate to one of $n-1$ other possible values picked uniformly at random (the chosen coordinate cannot stay the same). Starting from the 0 vector, the number of nonzero coordinates evolves as a birth and death chain on $\{0, \ldots, r\}$, with

$$p_x = \frac{(r-x)}{r}, \qquad q_x = \frac{x}{r(n-1)}, \qquad r_x = \frac{x(n-2)}{r(n-1)}.$$

The stationary distribution is

$$\nu_{n,r}(x) = \binom{r}{x}(n-1)^x n^{-r}.$$

This chain has eigenvalues (see, e.g., [19], Section 5, for an elementary argument)

$$\lambda_{n,r,i} = \frac{in}{r(n-1)}.$$

Hence,

$$\lambda_{n,r} = \frac{n}{r(n-1)},$$

$$t_{n,r} = \frac{r(n-1)}{n} \sum_1^r \frac{1}{i},$$

$$\theta_2(n,r) = \lambda_{n,r}^{-2} \sigma_{n,r}^2 = \sum_1^r \frac{1}{i^2}.$$

Fix a sequence $(n_\ell, r_\ell)$ with $r_\ell \to \infty$. Then by Theorem 5.1, the associated continuous-time chain has a $(s_\ell, \xi_\ell)$-separation cut-off with

$$s_\ell = (1 - 1/n_\ell) r_\ell \log r_\ell, \qquad \xi_\ell = (1 - 1/n_\ell) r_\ell.$$

The shape is given by the Gumbel distribution.

As a variation, consider the birth and death chain on $\{0, \ldots, r\}$ with

(7.1) $\quad p_x = \dfrac{r-x}{r}, \qquad q_x = \dfrac{x}{r}\theta, \qquad r_x = \dfrac{x}{r}(1-\theta), \qquad \theta \in (0,1).$

This has stationary distribution

$$\nu_{\theta,r}(x) = \binom{r}{x}\theta^{r-x}(1+\theta)^{-r}.$$

This is the projection (under the natural action of the symmetric group $S_r$) of the probability measure

$$\overline{\nu}_{\theta,r}(\overline{x}) = \frac{\theta^{|\overline{x}|}}{(1+\theta)^r}$$

on the hypercube $\{0,1\}^r$ with $\overline{x} = (x_i)_1^r$, $x_i \in \{0,1\}$ and $|\overline{x}| = \sum x_i$. The birth and death chain on $\{0, \ldots, r\}$ with rates (7.1) is the projection of a chain $K$ on the hypercube with $K(\overline{x}, \overline{x}) = |\overline{x}|(1-\theta)/r$,

$$K(\overline{x}, \overline{y}) = \begin{cases} 1/r, & \text{if } |\overline{x}| = |\overline{y}| + 1, \\ \theta/r, & \text{if } |\overline{x}| = |\overline{y}| - 1 \end{cases}$$

and $K(\overline{x}, \overline{y}) = 0$ otherwise (see, e.g., [17], Section 3). The eigenvalues are

$$\lambda_{\theta,r,i} = \frac{i}{r}(1+\theta).$$

Hence,

$$\lambda_{\theta,r} = \frac{1+\theta}{r}, \qquad t_{\theta,r} = \frac{r}{1+\theta}\sum_1^r \frac{1}{i}, \qquad \theta_2(n,r) = \lambda_{\theta,r}^{-2}\sigma_{\theta,r}^2 = \sum_1^r \frac{1}{i^2}.$$

Consider a sequence $(\theta_\ell, r_\ell)$ with $r_\ell$ tending to infinity. By Theorem 5.1, the associated continuous-time chain has a separation cut-off of type

$$((1+\theta_\ell)^{-1}r_\ell \log r_\ell, r_\ell)$$

[Theorem 5.1 gives a window of size $(1+\theta_\ell)^{-1}r_\ell$, but this is essentially equivalent since $\theta_\ell \in (0,1)$]. This chain is studied in [17] and [18], Section 5 (as a certain Metropolis chain on the hypercube, in discrete time). There it is proved that the chain has an $L^2$ cut-off of type

$$((1+\theta_\ell)^{-1}r_\ell \log \sqrt{r_\ell/\theta_\ell}, r_\ell).$$

The reference [13] also proves that for a fixed $\theta$, there is a

$$((1+\theta)^{-1}r \log \sqrt{r/\theta}, r)$$

total variation cut-off as $r$ tends to infinity. Note, however, that this last result cannot hold true if $\theta$ is allowed to vary and tend to 0. In general, for a sequence $(\theta_\ell, r_\ell)$ with $r_\ell$ tending to infinity, there is a total variation cut-off of type $(\rho_\ell, r_\ell)$ with

$$\rho_\ell = (1+\theta_\ell)^{-1}r_\ell \min\{\log r_\ell, \log \sqrt{r_\ell/\theta_\ell}\}.$$

As for the Bernoulli–Laplace models, we only need to prove a total variation lower bound matching the upper bound given by the separation and $L^2$ results. Such a total variation lower bound is easily derived using the data and method of [13], page 179.

*Distance regular graphs.* A finite graph is distance-transitive if the automorphism group of that graph acts transitively on the set of vertex pairs $(x,y)$ with $d(x,y) = k$, for any $k$. Distance-regular graphs generalize this notion without requiring a group action (see, e.g., [5, 7, 9]). Let $(V,E)$ be a connected graph with vertex set $V$ and symmetric edge set $E \subset V \times V$. Let $\rho$ be the graph distance and $m$ be the diameter of $(V,E)$. A graph is *regular* if the number of vertices at distance 1 from $x$ is independent of $x \in V$. A graph is *distance-regular* if for any $i,j \in \{0,\ldots,m\}$, the number of vertices at distance $i$ from $x$ and $j$ from $y$ depends only on the distance $\rho(x,y)$ between $x$ and $y$, $x,y \in V$.

One well-known basic result is that simple random walk on a distance-regular graph can be studied by collapsing to a birth and death chain on $\{0,\ldots,m\}$ started at 0; see [7], Section 3, where more general walks on distance-regular graphs that can be treated in a similar fashion are discussed. In this collapse, the set of all eigenvalues (without multiplicities) is conserved and separation, total variation or $L^p$ distance to the stationary measure for any fixed $p$ are conserved. Thus Theorem 5.1 yields the following result:

THEOREM 7.1. *Let $\mathcal{G}_n$ be a family of distance-regular graphs with diameter $m_n$ tending to infinity as $n$ tends to infinity. On the vertex set $V_n$ of $\mathcal{G}_n$, consider the continuous-time simple random walk and let $\gamma^t_{n,x_n}$ be the law of that process started at a fixed arbitrary point $x_n \in V_n$. Let $\nu_n$ be the uniform probability measure on $(V_n)$. Let $\lambda_n$ be the associated spectral gap and for any fixed $\varepsilon \in (0,1)$, set*

$$\tau_n = \inf\{t > 0 : \mathrm{sep}(\gamma^t_{n,x_n},\nu_n) \leq \varepsilon\}.$$

*Then the family $(V_n, \gamma^t_{n,x_n}, \nu_n)$ has a separation cut-off if and only if $\lambda_n \tau_n \to \infty$.*

REMARKS. 1. The separation $\mathrm{sep}(\gamma^t_{n,x_n},\nu_n)$ does not depend on the starting point $x_n$, so this result can be read as a max-separation result.

2. Note that the family $\mathcal{G}_n$ does not need to be "natural" in any way. It can mix elements from the various natural families described below.

Theorem 5.1 leads to the computation of the cut-off time, when it exists, in terms of the spectrum. It is conjectured by experts that distance-regular graphs have been classified. They include the following examples:

(i) The finite circle $\mathbb{Z}/n\mathbb{Z}$ with an edge from $x$ to $y$ if and only if $|x-y| = 1$. This family has no cut-off.

(ii) Hamming distance graphs such as the hypercube. These have been discussed above.

(iii) The natural graph on $r$-sets of an $n$-set with an edge from $x$ to $y$ if $\#(x \cap y) = k-1$. This is equivalent to the Bernoulli–Laplace models discussed above.

(iv) $q$-Families: These are described in some detail in [7] with data that is useful for our purpose. These families are all related to certain types of vector-subspaces of a finite-dimensional vector space over a finite field $\mathbb{F}_q$ ($q$ a prime power), hence the name. The simplest example is the set of all $m$-dimensional vector subspaces of an $n$-dimensional vector space discussed below.

Regarding the known distance-regular graphs from the $q$-Families ($q$-DRG for short) listed in [7], we can state the following theorem:

THEOREM 7.2. *Referring to the setting of Theorem 7.1, assume that the family $\mathcal{G}_n$ is made of some of the known $q$-DRGs listed in [7]. Then there is a separation cut-off if and only if the diameter $m_n$ of $\mathcal{G}_n$ tends to infinity. Moreover, if $m_n$ tends to infinity, then there is a $(m_n, \sqrt{m_n})$ separation cut-off with a nondegenerate normal shape.*

This immediately follows from the data reviewed in [7] and from Theorem 5.1. Instead of going into the details, we will illustrate the result with the simplest case of a natural $q$-DRG family.

*Fixed-dimension subspaces of a $q$ vector space.* Let $q$ be a prime power and $\mathbb{F}_q$ be a finite field of order $q$. Let $E_n$ be an $n$-dimensional vector space over $\mathbb{F}_q$. For $m \le n/2$, let $V_{q,n,m}$ be the set of all $m$-dimensional vector subspaces of $E_n$. The finite set $V_{q,n,m}$ is equipped with the distance

$$d(x,y) = m - \dim(x \cap y)$$

and the graph structure according to which $(x,y)$ is an edge if and only if $d(x,y) = 1$. The induced graph distance is the distance $d$ and this graph has diameter $m \le n/2$. The action of $GL_n(\mathbb{F}_q)$ on vector subspaces shows that this is a distance-transitive graph, hence a distance-regular graph. We consider the simple random walk on this graph (started at an arbitrary subspace) in continuous time. Let $\gamma_{q,n,m}^t$ be its law at time $t$ and $\nu_{q,n,m}$ its (uniform) stationary measure. Details (dealing with the discrete-time version) can be found in [7, 13]. As explained in [7, 13] and briefly above, this process can be studied through a birth and death chain on $\{0, \ldots, m\}$ which is simply the associated "distance process." Known computations involving relevant families of orthogonal polynomials give the eigenvalues $\lambda_{q,n,m,i}$ of this chain as

$$\lambda_{q,n,m,i} = \frac{(1-q^{-i})(1-q^{i-n-1})}{(1-q^{m-n})(1-q^{-m})}, \qquad 0 \le i \le m.$$

Hence the smallest nonzero eigenvalue is

$$\frac{(1-q^{-1})(1-q^{-n})}{(1-q^{m-n})(1-q^{-m})}$$

and when $m$ (hence also $n$) tends to infinity, we have

$$\sum_{i=1}^{m} \lambda_{q,n,m,i}^{-1} = m + O(1), \qquad \sum_{i=1}^{m} \lambda_{q,n,m,i}^{-2} = m + O(1).$$

If we now choose an arbitrary sequence $(q_\ell, n_\ell, m_\ell)$ and consider the family $(\Omega_\ell, \gamma_\ell^t, \nu_\ell)$ where

$$\Omega_\ell = V_{q_\ell, n_\ell, m_\ell}, \qquad \gamma_\ell^t = \gamma_{q_\ell, n_\ell, m_\ell}^t, \qquad \nu_\ell = \nu_{q_\ell, n_\ell, m_\ell},$$

then Theorem 5.1 shows that this family has a separation cut-off if and only if $m_\ell$ tends to infinity. Assuming that $m_\ell$ tends to infinity, Theorem 6.1 shows that there is an $(m_\ell, \sqrt{m_\ell})$-separation cut-off and that the window size $\sqrt{m_\ell}$ is optimal.

The references [7, 13], when translated into continuous time (in the present case, there are significant differences between discrete and continuous time), give a total variation cut-off of type $(m_\ell, \sqrt{m_\ell})$ and an $L^2$ cut-off of type $(s_\ell, 1)$ with

$$s_\ell = \tfrac{1}{2} m_\ell (n_\ell - m_\ell) \log q_\ell.$$

**Acknowledgment.** This work began at the American Institute of Mathematics in Palo Alto in December 2004 during the ARCC workshop "Sharp Thresholds for Mixing Times." We thank the organizers, Amir Dembo, Yuval Peres, David Revelle and the other participants of the workshop for their helpful comments.

DEPARTMENT OF STATISTICS
STANFORD UNIVERSITY
SEQUOIA HALL
390 SERRA MALL
STANFORD, CALIFORNIA 94305-4065
USA
E-MAIL: diaconis@math.stanford.edu

DEPARTMENT OF MATHEMATICS
CORNELL UNIVERSITY
MALOTT HALL
ITHACA, NEW YORK 14853-4201
USA
E-MAIL: lsc@math.cornell.edu